\documentclass[10pt]{amsart}
\usepackage{epic}

\theoremstyle{plain}
\newtheorem{thm}{Theorem}[section]
\newtheorem{lem}[thm]{Lemma}

\newtheorem{qn}{QUESTION}

\newtheorem{clm}[thm]{Claim}

\theoremstyle{definition}

\newcommand{\Q}{{\bf{Q}}}

\newcommand{\Hecke}{{\mathcal{H}}}
\newcommand{\Sym}{{\mathfrak{S}}}

\newcommand{\id}{{\mathrm{id}}}

\newcommand{\tr}{\operatorname{tr}}
\newcommand{\rad}{\operatorname{rad}}

\newcommand{\co}{\colon\thinspace}

\newcommand{\lplus}
{\raisebox{-23pt}
  {
  \begin{picture}(50,50)
  \put(0,0){\vector(1,1){50}}
  \drawline(50,0)(27,23)
  \put(23,27){\vector(-1,1){23}}
  \end{picture}
  }
}
\newcommand{\lminus}
{\raisebox{-23pt}
  {
  \begin{picture}(50,50)
  \put(50,0){\vector(-1,1){50}}
  \drawline(0,0)(23,23)
  \put(27,27){\vector(1,1){23}}
  \end{picture}
  }
}
\newcommand{\lzero}
{\raisebox{-23pt}
  {
  \begin{picture}(50,50)
  \qbezier(0,0)(20,20)(20,25)
  \qbezier(20,25)(20,30)(0,50)
  \qbezier(50,0)(30,20)(30,25)
  \qbezier(30,25)(30,30)(50,50)
  \put(0,50){\vector(-1,1){0}}
  \put(50,50){\vector(1,1){0}}
  \end{picture}
  }
}
\newcommand{\aA}
{\raisebox{-28pt}
  {
  \begin{picture}(40,60)
  \put(0,0){\line(1,1){18}}
  \qbezier(22,22)(30,30)(22,38)
  \put(18,42){\vector(-1,1){18}}

  \put(40,0){\line(-1,1){20}}
  \qbezier(20,20)(10,30)(20,40)
  \put(20,40){\vector(1,1){20}}

  \end{picture}
  }
}
\newcommand{\oneone}
{\raisebox{-28pt}
  {
  \begin{picture}(40,60)
  \qbezier(0,0)(30,30)(0,60)
  \put(0,60){\vector(-1,1){0}}

  \qbezier(40,0)(10,30)(40,60)
  \put(40,60){\vector(1,1){0}}
  \end{picture}
  }
}
\newcommand{\aba}
{\raisebox{-29pt}
  {
  \begin{picture}(40,62)
  \qbezier(0,0)(0,10)(10,20)
  \drawline(10,20)(30,40)
  \qbezier(30,40)(40,50)(40,60)
  \put(40,60){\vector(0,1){2}}

  \qbezier(20,0)(20,10)(12,18)
  \qbezier(8,22)(0,30)(8,38)
  \drawline(8,38)(12,42)
  \qbezier(12,42)(20,50)(20,60)
  \put(20,60){\vector(0,1){2}}

  \qbezier(40,0)(40,10)(30,20)
  \drawline(30,20)(22,28)
  \drawline(18,32)(12,38)
  \qbezier(8,42)(0,50)(0,60)
  \put(0,60){\vector(0,1){2}}
  \end{picture}
  }
}

\newcommand{\bab}
{\raisebox{-29pt}
  {
  \begin{picture}(40,62)
  \qbezier(0,0)(0,10)(10,20)
  \drawline(10,20)(30,40)
  \qbezier(30,40)(40,50)(40,60)
  \put(40,60){\vector(0,1){2}}

  \qbezier(20,0)(20,10)(28,18)
  \drawline(28,18)(32,22)
  \qbezier(32,22)(40,30)(32,38)
  \qbezier(28,42)(20,50)(20,60)
  \put(20,60){\vector(0,1){2}}

  \qbezier(40,0)(40,10)(32,18)
  \drawline(28,22)(22,28)
  \drawline(18,32)(10,40)
  \qbezier(10,40)(0,50)(0,60)
  \put(0,60){\vector(0,1){2}}
  \end{picture}
  }
}

\newcommand{\aCab}
{\raisebox{-48pt}
  {
  \begin{picture}(60,100)
  \qbezier(0,0)(0,30)(8,38)
  \drawline(8,38)(12,42)
  \qbezier(12,42)(20,50)(20,60)
  \qbezier(20,60)(20,70)(12,78)
  \qbezier(8,82)(0,90)(0,100)

  \qbezier(20,0)(20,10)(28,18)
  \drawline(28,18)(32,22)
  \qbezier(32,22)(40,30)(40,40)
  \qbezier(40,40)(40,50)(48,58)
  \qbezier(52,62)(60,70)(60,100)

  \qbezier(40,0)(40,10)(32,18)
  \drawline(28,22)(12,38)
  \qbezier(8,42)(0,50)(0,60)
  \qbezier(0,60)(0,70)(8,78)
  \drawline(8,78)(12,82)
  \qbezier(12,82)(20,90)(20,100)

  \qbezier(60,0)(60,50)(52,58)
  \drawline(52,58)(48,62)
  \qbezier(48,62)(40,70)(40,100)
  
  \end{picture}
  }
}

\newcommand{\one}
{\raisebox{-49pt}
  {
  \begin{picture}(70,102)
  \put(0,0){\vector(0,1){102}}
  \put(20,0){\vector(0,1){102}}
  \put(40,0){\vector(0,1){102}}
  \put(60,0){\vector(0,1){102}}
  \put(60,50){\circle*{2}}
  \put(62,45){$x_0$}
  \end{picture}
  }
}

\newcommand{\pullleft}
{\raisebox{-49pt}
  {
  \begin{picture}(80,102)
  \put(20,0){\line(0,1){38}}
  \put(20,42){\line(0,1){16}}
  \put(20,62){\vector(0,1){40}}
  \put(40,0){\line(0,1){28}}
  \put(40,32){\line(0,1){36}}
  \put(40,72){\vector(0,1){30}}
  \put(60,0){\line(0,1){18}}
  \put(60,22){\line(0,1){56}}
  \put(60,82){\vector(0,1){20}}
  \qbezier(80,0)(80,10)(60,20)
  \put(60,20){\line(-2,1){40}}
  \qbezier(20,40)(10,45)(10,50)
  \put(10,50){\circle*{2}}
  \put(0,45){$x_0$}
  \qbezier(10,50)(10,55)(20,60)
  \put(20,60){\line(2,1){40}}
  \qbezier(60,80)(80,90)(80,100)
  \put(80,100){\vector(0,1){2}}
  \end{picture}
  }
}

\title{Braid groups and Iwahori-Hecke algebras}
\thanks{Partially supported by NSF grant DMS 0307235 and the Sloan Foundation.}
\author{Stephen Bigelow}
\address{Department of Mathematics,
         University of California at Santa Barbara,
         California 93106, USA}
\email{bigelow@math.ucsb.edu}
\date{April 2005}

\begin{document}

\begin{abstract}
The braid group $B_n$ is
the mapping class group of an $n$-times punctured disk.
The Iwahori-Hecke algebra $\Hecke_n$ is a quotient of
the braid group algebra of $B_n$
by a quadratic relation in the standard generators.
We discuss how to use $\Hecke_n$
to define the Jones polynomial of a knot or link.
We also summarize
the classification of the irreducible representations of $\Hecke_n$.
We conclude with some directions for future research
that would apply mapping class group techniques
to questions related to $\Hecke_n$.
\end{abstract}

\maketitle

\section{Introduction}

The braid group $B_n$
is the mapping class group of an $n$-times punctured disk.
It can also be defined using
certain kinds of arrangements of strings in space,
or certain kinds of diagrams in the plane.
Our main interest in the braid group $B_n$
will be in relation to the Iwahori-Hecke algebra $\Hecke_n$,
which is a certain quotient of the group algebra of $B_n$.
The exact definition will be given in Section \ref{sec:hecke}.

The Iwahori-Hecke algebra plays
an important role in representation theory.
It first came to the widespread attention of topologists
when Jones used it to define the knot invariant
now called the Jones polynomial \cite{vJ85}.
This came as a huge surprise,
since it brought together two subjects
that were previously unrelated.
It has given knot theorists a host of new knot invariants,
and intriguing connections to other areas of mathematics to explore.
It has also helped to promote
the use of pictures and topological thinking in
representation theory.

As far as I know,
no major results related to the Iwahori-Hecke algebra
have yet been proved
using the fact that $B_n$ is a mapping class group.
I think the time is ripe for such a result.
Unfortunately this paper will 
use the diagrammatic definition of the braid group almost exclusively.
I hope it will at least help to provide a basic grounding
for someone who wants to pursue
the connection to mapping class groups in the future.

The outline of this paper is as follows.
In Sections \ref{sec:braid} and \ref{sec:hecke}
we introduce the braid group and Iwahori-Hecke algebra.
In Sections \ref{sec:basis} and \ref{sec:trace}
we give a basis for the Iwahori-Hecke algebra
and for its module of trace functions.
In Section \ref{sec:jones}
we explain how one such trace function leads to the definition of
the Jones polynomial of a knot or a link.
In Section \ref{sec:rep} we briefly summarize
the work of Dipper and James \cite{DJ86}
classifying the irreducible representations of the Iwahori-Hecke algebra.
In Section \ref{sec:future}
we conclude with some speculation on possible directions for future research.
Open problems will be scattered throughout the paper.

\section{The braid group}
\label{sec:braid}

Like most important mathematical objects,
the braid group $B_n$ has several equivalent definitions.
Of greatest relevance to this volume
is its definition as a mapping class group.
Let $D$ be a closed disk,
let $p_1,\dots,p_n$ be distinct points in the interior of $D$,
and let $D_n = D \setminus \{p_1,\dots,p_n\}$.
The braid group $B_n$ is the mapping class group of $D_n$.
Thus a braid is the equivalence class of a homeomorphism from $D_n$ to itself
that acts as the identity of the boundary of the disk.

Artin's original definition of $B_n$
was in terms of {\em geometric braids}.
A geometric braid is a disjoint union of
$n$ edges, called {\em strands}, in $D \times I$,
where $I$ is the interval $[0,1]$.
The set of endpoints of the strands is required to be
$\{p_1,\dots,p_n\} \times \{0,1\}$,
and each strand is required to intersect each disk cross-section exactly once.
Two geometric braids are said to be equivalent
if it is possible to deform one to the other
through a continuous family of geometric braids.
The elements of $B_n$ are equivalence classes of geometric braids.

We will need some terminology to refer to directions in a geometric braid.
Take $D$ to be the unit disk centered at $0$ in the complex plane.
Take $p_1,\dots,p_n$ to be real numbers
with $-1 < p_1 < \dots < p_n < 1$.
The {\em top} and {\em bottom} of the braid are
$D \times \{1\}$ and $D \times \{0\}$, respectively.
In a disk cross-section,
the {\em left} and {\em right}
are the directions of decreasing and increasing real part, respectively,
while the {\em front} and {\em back}
are the directions of decreasing and increasing imaginary part, respectively.

Multiplication in $B_n$ is defined as follows.
If $a$ and $b$ are geometric braids with $n$ strands
then the product $ab$ is obtained
by stacking $a$ on top of $b$
and then rescaling vertically to the correct height.
This can be shown to give a well-defined product of equivalence classes,
and to satisfy the axioms of a group.

\begin{figure}
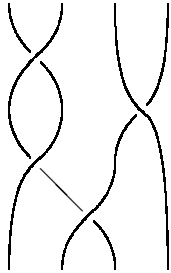

$$\aCab$$
\caption{A braid with four strands.}
\label{fig:fourbraid}
\end{figure}
A geometric braid can be drawn in the plane using
a projection from $D \times I$ to $[-1,1] \times I$.
An example is shown in Figure \ref{fig:fourbraid}.
The projection map is given by $(x+iy,t) \mapsto (x,t)$.
Note that this sends each strand to an embedded edge.
We also require that the braid be in {\em general position}
in the sense that
the images of the strands intersect each other transversely,
with only two edges meeting at each point of intersection.
The points of intersection are called {\em crossings}.
At each crossing,
we record which of the two strands passed in front of the other
at the corresponding disk cross-section of the geometric braid.
This is usually represented pictorially
by a small break in the segment that passes behind.
The image of a geometric braid under a projection in general position,
together with this crossing information, is called a {\em braid diagram}.

Let us fix some terminology related to braid diagrams.
The directions left, right, top, and bottom
are the images of these same directions in the geometric braid,
so for example the point $(-1,1)$ is the top left of the braid diagram.
We say a strand makes an {\em overcrossing} or an {\em undercrossing}
when it passes respectively
in front of or behind another strand at a crossing.
A crossing is called {\em positive}
if the strand making the overcrossing
goes from the bottom left to the top right of the crossing,
otherwise it is called {\em negative}.
The endpoints of strands are called {\em nodes}.

\begin{figure}
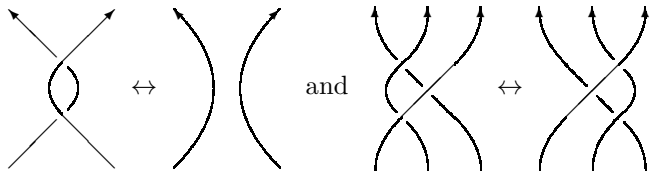

$$\aA \leftrightarrow \oneone \ {\mathrm{\ and\ }} \ \aba \leftrightarrow \bab$$
\caption{Reidemeister moves of types two and three.}
\label{fig:reidemeister}
\end{figure}
Two braid diagrams represent the same braid
if and only if they are related by an isotopy of the plane
and a sequence of Reidemeister moves of types two and three.
These are moves in which the diagram remains unchanged
except in a small disk,
where it changes as shown in Figure \ref{fig:reidemeister}.
(There is also a Reidemeister move of type one,
which is relevant to knots but not to braids.)

For $i=1,\dots,n-1$,
let $\sigma_i$ be the braid diagram with one crossing,
which is a positive crossing between strands $i$ and $i+1$.
The braid group $B_n$ is
generated by $\sigma_1,\dots,\sigma_{n-1}$,
with defining relations
\begin{itemize}
\item $\sigma_i \sigma_j = \sigma_j \sigma_i$ if $|i-j| > 1$,
\item $\sigma_i \sigma_j \sigma_i = \sigma_j \sigma_i \sigma_j$ if $|i-j| = 1$.
\end{itemize}

There is an imprecise but vivid physical description of
the correspondence between a geometric braid
and a mapping class of the $n$-times punctured disk $D_n$.
Imagine a braid made of inflexible wires
and a disk made of flexible rubber.
Press the disk onto the top of the braid,
puncturing the disk at $n$ points in its interior.
Now hold the disk by its boundary and push it down.
As the wires of the braid twist around each other,
the punctures of the disk will twist around
and the rubber will be stretched and distorted to accommodate this.
The mapping class corresponding to the geometric braid
is represented by the function taking each point on $D_n$
to its image in $D_n$ after the disk has been pushed
all the way to the bottom of the braid.
(With our conventions,
this description gives the group of mapping classes acting on the right.)

See \cite{jB74}, or \cite{BB05}, for proofs that
these and other definitions of $B_n$ are all equivalent.
This paper will primarily use
the definition of a braid group as a braid diagram.
This is in some sense the least elegant choice
since it involves an arbitrary projection
and a loss of the true three-dimensional character of the geometric braid.
The main goal is to provide an introduction
that may inspire someone to apply mapping class group techniques
to problems that have previously been studied
algebraically and combinatorially.

\section{The Iwahori-Hecke algebra}
\label{sec:hecke}

Let $n$ be a positive integer
and let $q_1$ and $q_2$ be units in a domain $R$.
The {\em Iwahori-Hecke algebra} $\Hecke_n(q_1,q_2)$,
or simply $\Hecke_n$,
is the associative $R$-algebra given by
generators $T_1,\dots,T_{n-1}$ and relations
\begin{itemize}
\item $T_i T_j = T_j T_i$ if $|i-j| > 1$,
\item $T_i T_j T_i = T_j T_i T_j$ if $|i-j| = 1$,
\item $(T_i - q_1)(T_i - q_2) = 0$.
\end{itemize}

The usual definition of the Iwahori-Hecke algebra
uses only one parameter $q$.
It corresponds to $\Hecke_n(-1,q)$,
or in some texts to $\Hecke_n(1,-q)$.
There is no loss of generality because there is an isomorphism
from $\Hecke_n(q_1,q_2)$ to $\Hecke_n(-1,-q_2/q_1)$
given by $T_i \mapsto -q_1 T_i$.
It will be convenient for us to keep two parameters.

We now explore some of the basic properties of
the Iwahori-Hecke algebra.
Since $q_1$ and $q_2$ are units of $R$,
the generators $T_i$ are units of $\Hecke_n$, with
$$T_i^{-1} = (T_i - q_1 - q_2)/(q_1 q_2).$$
Thus there is a well-defined homomorphism
from $B_n$ to the group of units in $\Hecke_n$ given by
$$\sigma_i \mapsto T_i.$$
The following is a major open question.

\begin{qn}
\label{qn:inj}
If $R = \Q(q_1,q_2)$,
is the above map from $B_n$ to $\Hecke_n(q_1,q_2)$ injective?
\end{qn}

For $n=3$, the answer is yes.
For $n=4$, the answer is yes if and only if
the Burau representation of $B_4$ is injective, or {\em faithful}.
The Burau representation is
one of the irreducible summands of the Iwahori-Hecke algebra
over $\Q(q_1,q_2)$.
By a result of Long \cite{dL86},
the map from $B_n$ to $\Hecke_n$ is injective if and only if
at least one of these irreducible summands is faithful.
For $n=4$,
they are all easily shown to be unfaithful
except for the Burau representation, which remains unknown.
For $n \ge 5$,
the Burau representation is unfaithful \cite{sB99},
but there are other summands whose status remains unknown.

One can also ask Question \ref{qn:inj} for other
choices of ring $R$ and parameters $q_1$ and $q_2$.
If the map from $B_n$ to $\Hecke_n$ is injective
for any such choice
then it is injective when $R = \Q(q_1,q_2)$.
A non-trivial case when the map is not injective is
when $n=4$ and $R = k[q_1^{\pm 1},q_2^{\pm 1}]$,
where $k$ is a field of characteristic $2$ \cite{CL97} or $3$ \cite{CL98}.
Another is when $n=4$, $R = \Q$
and $q_2 / q_1 = -2$ \cite{sB00}.

\begin{figure}
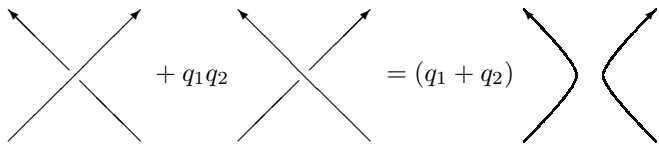

$$\lplus + q_1 q_2 \lminus = (q_1 + q_2) \lzero$$
\caption{The skein relation}
\label{fig:skein}
\end{figure}
Using the map from $B_n$ to $\Hecke_n$,
we can represent any element of $\Hecke_n$
by a linear combination of braid diagrams.
The quadratic relation is equivalent to
the {\em skein relation} shown in Figure \ref{fig:skein}.
Here, an instance of the skein relation
is a relation involving three diagrams
that are identical except inside a small disk
where they are as shown in the figure.

One motivation for studying the Iwahori-Hecke algebra
is its connection with the representation theory of the braid groups.
The representations of $\Hecke_n$
are precisely those representations of $B_n$
for which the image of the generators satisfy a quadratic relation.
The study of these representations
led Jones to the discovery of his knot invariant,
which we define in Section \ref{sec:jones}

Another reason for interest is the connection between
the Iwahori-Hecke algebra and the symmetric group.
There is an isomorphism from $\Hecke_n(1,-1)$ to the group algebra $R \Sym_n$
taking $T_i$ to the transposition $(i,i+1)$.
Thus $\Hecke_n(q_1,q_2)$ can be thought of as
a {\em deformation} of $R \Sym_n$.
The Iwahori-Hecke algebra plays a role in
the representation theory of the general linear group over a finite field
that is analogous to the role of
the symmetric group in the representation theory of
the general linear group over the real numbers.
See for example \cite{rD85}.

This process of
realizing a classical algebraic object as the case $q=1$
in a family of algebraic objects parametrized by $q$
is part of a large circle of ideas called {\em quantum mathematics},
or {\em $q$-mathematics}.
The exact nature and significance of any connection to quantum mechanics
not clear at present.
One example is \cite{jB01},
in which Barrett uses quantum mathematics to analyze quantum gravity
in a universe with no matter and three space-time dimensions.

\section{A basis}
\label{sec:basis}

The aim of this section is to show
that $\Hecke_n$ is a free $R$-module of rank $n!$,
and to give an explicit basis.

Let $\phi \co B_n \to \Sym_n$ be the map such that
$\phi(\sigma_i)$ is the transposition $(i,i+1)$.
Thus in any braid $b$,
the strand with lower endpoint at node number $i$
has upper endpoint at node number $\phi(b)(i)$.

For $w \in \Sym_n$,
let $T_w$ be a braid diagram with the minimal number of crossings
such that every crossing is positive and $\phi(T_w) = w$.
Such a braid can be thought of as ``layered'' in the following sense.
In the front layer is a strand connecting
node $1$ at the bottom to node $w(1)$ at the top.
Behind that is a strand connecting node $2$ at the bottom
to node $w(2)$ at the top.
This continues until the back layer, in which
a strand connects node $n$ at the bottom to node $w(n)$ at the top.
From this description it is clear that
our definition of $T_w$ specifies a unique braid in $B_n$.
By abuse of notation,
let $T_w$ denote the image of this braid in $\Hecke_n$.
For example, if $w$ is a transposition $(i,i+1)$
then $T_w$ is the generator $T_i$.

\begin{thm}
\label{thm:heckebasis}
The set of $T_w$ for $w \in \Sym_n$ forms a basis for $\Hecke_n$.
\end{thm}

To prove this, we first describe an algorithm
that will input a linear combination of braid diagrams
and output a linear combination of basis elements $T_w$
that represents the same element of $\Hecke_n$.
By linearity,
it suffices to describe
how to apply the algorithm to a single braid diagram $v$.

A crossing in $v$ will be called {\em bad} 
if the strand that makes the overcrossing
is the one whose lower endpoint is farther to the right.
If $v$ has no bad crossings, stop here.

Suppose $v$ has at least one bad crossing.
Let the {\em worst} crossing be a bad crossing
whose undercrossing strand has lower endpoint farthest to the left.
If there is more than one such bad crossing,
let the worst be the one that is closest to the bottom of the diagram.

Use the skein relation to rewrite $v$
as a linear combination of $v'$ and $v_0$,
where $v'$ is the result of changing the sign of the worst crossing
and $v_0$ is the result of removing it.
Now recursively apply this procedure to $v'$ and $v_0$.

Note that any bad crossings in $v'$ and $v_0$
are ``better'' than the worst crossing of $v$ in the sense that either
the lower endpoint of their undercrossing strand is farther to the right
or they have the same undercrossing strand
and are closer to the top of the diagram.
Thus the above algorithm must eventually terminate
with a linear combination of diagrams that have no bad crossings.
Any such diagram must equal $T_w$ for some $w \in \Sym_n$.

This algorithm shows that the $T_w$ span $\Hecke_n$.
It remains to show that they are linearly independent.
Note that if the algorithm is given as input
a linear combination of diagrams of the form $T_w$,
then its output will be the same linear combination.
Thus it suffices to show that
the output of the algorithm does not depend on the initial choice of
linear combination of braid diagrams
to represent a given element of $\Hecke_n$.
We prove this in three claims,
which show that the output of the algorithm is invariant under
the skein relation and Reidemeister moves of types two and three.

\begin{clm}
\label{clm:skein}
Suppose $v_+$, $v_-$ and $v_0$ are three braid diagrams
that are identical except in a small disk
where $v_+$ has a positive crossing,
$v_-$ has a negative crossing, and
$v_0$ has no crossing.
Then the algorithm gives the same output
for both sides of the skein relation
$v_+ + q_1q_2 v_- = (q_1 + q_2)v_0$.
\end{clm}

\begin{proof}
For exactly one of $v_+$ and $v_-$,
the crossing inside the small disk is a bad crossing.
For convenience assume it is $v_+$,
since it makes no difference to the argument.

Suppose the worst crossing for $v_+$
is the crossing in the small disk.
Applying the next step of the algorithm to $v_+$
results in a linear combination of $v_-$ and $v_0$
which, by design,
will exactly cancel the other two terms in the skein relation.

Now suppose the worst crossing for $v_+$ is not inside the small disk.
Then it must be the same as the worst crossing for $v_-$ and for $v_0$.
Thus the next step of the algorithm has the same effect on
$v_+$, $v_-$ and $v_0$.
The claim now follows by induction.
\end{proof}

\begin{clm}
\label{clm:rtwo}
If $u$ and $v$ are diagrams that differ by a Reidemeister move of type two
then the algorithm gives the same output for $u$ as for $v$.
\end{clm}

\begin{proof}
As in the proof of the previous claim,
we can reduce to the case where a worst crossing
lies inside the small disk affected by the Reidemeister move.
The claim now follows by
computing the result of applying the algorithm inside the small disk.
Alternatively,
observe that this computation amounts to checking
the case $n=2$ of Theorem \ref{thm:heckebasis},
which follows easily from the presentation of $\Hecke_2$.
\end{proof}

\begin{clm}
\label{clm:rthree}
Suppose $u$ and $v$ are diagrams
that differ by a Reidemeister move of type three.
Then the algorithm gives the same output for $u$ as for $v$.
\end{clm}

\begin{proof}
Once again,
one solution involves a brute force computation of the algorithm.
Here we describe a somewhat more comprehensible approach.

Label the three strands in the small disk in each of $u$ and $v$
the {\em front}, {\em back}, and {\em middle} strands,
where the front strand makes two overcrossings,
the back strand makes two undercrossings,
and the middle strand makes one overcrossing and one undercrossing.

Let $u'$ and $v'$ be the result of
changing the sign of the crossings
between the front and middle strands of $u$ and $v$ respectively.
Note that eliminating these crossings results in identical braid diagrams.
Thus by Claim \ref{clm:skein},
the algorithm gives the same output for $u$ as for $v$
if and only if it gives the same output for $u'$ as for $v'$.

Relabel the three strands in $u'$ and $v'$
so that once again the front strand makes two overcrossings,
the back strand makes two undercrossings,
and the middle strand makes one overcrossing and one undercrossing.
Now let $u''$ and $v''$ be the output of
changing the sign of the crossing
between the middle and back strands of $u'$ and $v'$ respectively.
Note that eliminating these crossings results in
braid diagrams that differ by Reidemeister moves of type two.
Thus by Claims \ref{clm:skein} and \ref{clm:rtwo},
the algorithm gives the same output for $u'$ as for $v'$
if and only if it gives the same output for $u''$ as for $v''$.

We can continue in this way,
alternately changing crossings between
front and middle, and middle and back strands.
We obtain six different versions of
the Reidemeister move of type three.
Each is obtained from the original by some crossing changes,
and corresponds to one of the six permutations of the roles of
front, middle, and back strands.

The algorithm gives the same output for $u$ as for $v$
if and only if it gives the same output
when the relevant disks in $u$ and $v$ are changed to represent
any one of the six versions of the Reidemeister move of type three.
Thus we can choose a version to suit our convenience.
In particular we can always choose
the front and back strands to be the ones with
lower endpoints farthest to the left and right respectively.
That way there will be no bad crossings inside the small disk,
and the algorithm will proceed identically
for the diagrams on either side of the move.

This completes the proof of the claim,
and hence of the theorem.
\end{proof}

\section{Trace functions}
\label{sec:trace}

A trace function on $\Hecke_n$
is a linear function $\tr \co \Hecke_n \to R$
such that $\tr(ab) = \tr(ba)$ for all $a, b \in \Hecke_n$.
Let $V$ be the quotient of $\Hecke_n$
by the vector subspace spanned by
elements of the form $ab-ba$ for $a, b \in \Hecke_n$.
Then the trace functions of $\Hecke_n$
correspond to the linear maps from $V$ to $R$.

The aim of this section is to find a basis for $V$,
and hence classify all trace functions of $\Hecke_n$.
This has been done by Turaev \cite{vT88}
and independently by Hoste and Kidwell \cite{HK90}.
They actually consider a larger algebra
in which the strands can have arbitrary orientations,
but the result is very similar.

We define a {\em closed $n$-braid} to be
a disjoint union of circles in $D \times S^1$
that intersects each disk cross-section at a total of $n$ points.
We say two closed braids are equivalent
if one can be deformed to the other
through a continuous family of closed braids.
The {\em closure} of a geometric braid in $D \times I$
is the result of identifying $D \times \{0\}$ to $D \times \{1\}$.
It is not difficult to show that
two braids have equivalent closures if and only if they are conjugate in $B_n$.

Define a {\em diagram} of a closed braid
to be a projection onto the annulus $I \times S^1$ in general position,
together with crossing information,
similar to the diagram of a braid.
Then $V$ is the vector space of
formal linear combination of closed $n$-braid diagrams
modulo the skein relation and Reidemeister moves of types two and three.
This is an example of a {\em skein algebra} of the annulus.

A {\em partition of $n$}
is a sequence $\lambda = (\lambda_1,\dots,\lambda_k)$
of integers such that $\lambda_1 \ge \dots \ge \lambda_k > 0$
and $\lambda_1 + \dots + \lambda_k = n$.
For any $m > 0$,
let $b_{(m)}$ be the braid
$\sigma_{m-1} \dots \sigma_2 \sigma_1$.
If $\lambda = (\lambda_1,\dots,\lambda_k)$ is a partition of $n$,
let $b_\lambda$ be the braid with diagram consisting of
a disjoint union of diagrams of the braids $b_{(\lambda_i)}$,
in order from left to right.
Let $v_\lambda$ be the closure of $b_\lambda$.

\begin{thm}
The set of $v_\lambda$ for partitions $\lambda$ of $n$ forms a basis for $V$.
\end{thm}

We start by defining an algorithm
similar to that of Theorem \ref{thm:heckebasis}.
There are some added complications because
a strand can circle around and cross itself,
and there is no ``bottom'' of the closed braid
to use as a starting point.
Therefore the first step is to
choose a basepoint $x_0$ on the diagram that is not a crossing point.
\begin{figure}
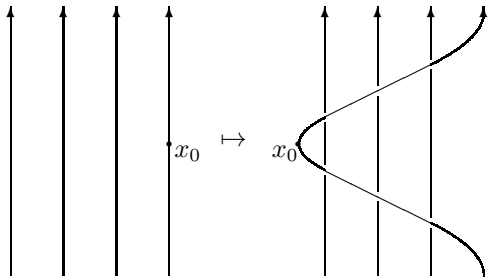

$$\one \ \mapsto \ \pullleft$$
\caption{Pulling the basepoint toward $\{0\} \times S^1$}
\label{fig:pullleft}
\end{figure}
Pull $x_0$ in front of the other strands so that $x_0$ becomes
the closest point to the boundary component $\{0\} \times S^1$,
as suggested by Figure \ref{fig:pullleft}.

Consider the oriented edge that begins at $x_0$
and proceeds in the positive direction around the annulus.
Call a crossing {\em bad} if this edge makes an undercrossing
on the first (or only) time it passes through that crossing.
If there is a bad crossing,
use the skein relation to eliminate
the first bad crossing the edge encounters.
This process will eventually terminate
with a linear combination of diagrams that have no bad crossings.

If there are no bad crossings then the loop through $x_0$
passes in front of every other loop in the closed braid.
Furthermore,
we can assume that its distance toward the front of the diagram
steadily decreases as it progresses in the positive direction from $x_0$
until just before it closes up again at $x_0$.
Since $x_0$ is the closest point to $\{0\} \times S^1$, this implies that
the loop through $x_0$ is isotopic to $v_{(m)}$ for some positive integer $m$.

Isotope this loop, keeping it in front of all other strands,
toward the boundary component $\{1\} \times S^1$,
until its projection is disjoint from that of all other loops.
Now ignore this loop
and repeat the above procedure
to the remainder of the closed braid diagram.
This process must eventually terminate with
a linear combination of closed braids of the form $v_\lambda$.

To show that the $v_\lambda$ are linearly independent,
it suffices to show that output is invariant under
the skein relation,
Reidemeister moves of types two and three,
and the choices of basepoint.
By induction on $n$ we can assume that
the output of the algorithm
does not depend on choices of basepoint made
after the first loop has been made disjoint from the other loops.
Thus the algorithm produces a unique output
given a diagram of a closed $n$-braid
together with a single choice of initial basepoint $x_0$.

\begin{clm}
The output of the algorithm is invariant under the skein relation,
and under any Reidemeister move
for which the basepoint does not lie in the disk affected by the move.
\end{clm}

\begin{proof}
The proofs of Claims \ref{clm:skein}, \ref{clm:rtwo} and \ref{clm:rthree}
go through unchanged.
\end{proof}

It remains only to prove the following.

\begin{clm}
For a given diagram of a closed $n$-braid,
the output of the algorithm does not depend on the choice of basepoint.
\end{clm}

\begin{proof}
First we show that the output of the algorithm is not affected
by moving the basepoint over an overcrossing.
Recall that the first step of the algorithm
is to pull the basepoint in front of the other strands
as in Figure \ref{fig:pullleft}.
If the basepoint is moved over an overcrossing,
the output of this first step will be altered by
a Reidemeister move of type two.
Furthermore, the basepoint does not lie in the disk
affected by this Reidemeister move.
Thus the output of the algorithm is unchanged.

Now fix a diagram $v$ of a closed $n$-braid.
By induction, assume that the claim is true
for any diagram with fewer crossings than $v$.
By the skein relation,
if $v'$ is the result of changing the sign
of one of the crossings of $v$
then the claim is true for $v$
if and only if it is true for $v'$.
Thus we are free to change the signs of the crossings in $v$
to suit our convenience.

By changing crossings and moving the basepoint past overcrossings
we can move the basepoint to any other point on the same loop in $v$.
Now suppose basepoints $x_0$ and $y_0$ 
lie on two distinct loops in $v$.
We can choose each to be
the closest point on its loop to the boundary component $\{0\} \times S^1$.
Assume, without loss of generality,
that $x_0$ is at least as close as to $\{0\} \times S^1$ as $y_0$ is.
By changing the signs of crossings,
we can assume that the loop through $y_0$ has no bad crossings.
Applying the algorithm using $y_0$ as the basepoint then has the affect of
isotoping the corresponding loop toward $\{1\} \times S^1$,
keeping it in front of all other strands.
This can be achieved by
a sequence of Reidemeister moves of types two and three.
No point in this loop is closer to $\{0\} \times S^1$ than $x_0$,
so $x_0$ does not lie in the disk affected by any of these Reidemeister moves.
Thus the algorithm will give the same output using the basepoint $x_0$
as it does using basepoint $y_0$.

This completes the proof of the claim,
and hence of the theorem.
\end{proof}

\section{The Jones polynomial}
\label{sec:jones}

The Jones polynomial is an invariant of knots and links,
first defined by Jones \cite{vJ85}.
Jones arrived at his definition
as an outgrowth of his work on operator algebras,
as opposed to knot theory.
To this day the topological meaning of his polynomial
seems somewhat mysterious,
and it has a very different flavor to classical knot invariants
such as the Alexander polynomial.

After the discovery of the Jones polynomial,
several people independently realized
that it could be generalized to a two-variable polynomial
now called the HOMFLY or HOMFLYPT polynomial.
The names are acronyms of the authors of
\cite{FYH85}, where the polynomial was defined,
and of \cite{PT88}, where related results were discovered independently.

The aim of this section is to show how to use the Iwahori-Hecke algebra
to define a polynomial invariant of knots and links
called the HOMFLY or HOMFLYPT polynomial.

Given a geometric braid $b$,
we can obtain a closed braid in the solid torus
by identifying the top and the bottom of $b$.
Now embed this solid torus into $S^3$
in a standard unknotted fashion.
The resulting knot or link in $S^3$
is called the {\em closure} of $b$.
It is a classical theorem of Alexander
that any knot or link in $S^3$ can be obtained in this way.

Let $B_\infty$ be the disjoint union of
the braid groups $B_n$ for $n \ge 1$.
For every $n \ge 1$, let
$$\iota \co B_n \to B_{n+1}$$
be the inclusion map that adds a single straight strand
to the right of any $n$-braid.
The {\em Markov moves} are as follows.
\begin{itemize}
\item $ab \leftrightarrow ba$,
\item $b \leftrightarrow \sigma_n \iota(b)$,
\item $b \leftrightarrow \sigma_n^{-1} \iota(b)$,
\end{itemize}
for any $a,b \in B_n$.

\begin{thm}[Markov's theorem]
Two braids have the same closure
if and only if they are connected by a sequence of Markov moves.
\end{thm}

An $R$-valued link invariant is thus equivalent to
a function from $B_\infty$ to $R$ that is invariant under the Markov moves.
We now look for such a function
that factors through the maps $B_n \to \Hecke_n$.

Let $\iota \co \Hecke_n \to \Hecke_{n+1}$ be
the inclusion map $T_i \mapsto T_i$.
A family of linear maps $\tr \co \Hecke_n \to R$ defines a link invariant
if and only if it satisfies the following.
\begin{itemize}
\item $\tr(ab) = \tr(ba)$,
\item $\tr(b) = \tr (T_n \iota(b))$,
\item $\tr(b) = \tr (T_n^{-1} \iota(b))$,
\end{itemize}
for every $n \ge 1$ and $a,b  \in \Hecke_n$.
We will call such a family of maps a {\em normalized Markov trace}.
The usual definition of Markov trace is slightly different,
but can easily be rescaled to satisfy the above conditions.

By the skein relation,
the third condition on a normalized Markov trace is equivalent to
\begin{equation}
\label{eq:addstrand}
(1+q_1 q_2)\tr(b) = (q_1 + q_2)\tr(\iota(b)),
\end{equation}
for every $n \ge 1$ and $b \in \Hecke_n$.
To obtain an interesting invariant,
we assume from now on that $q_1 + q_2$ is a unit of $R$.

Let $\lambda$ be a partition of $n$
and let $b_\lambda$ be the corresponding braid diagram
as defined in the previous section.
Let $k$ be the number of components of $b_\lambda$,
that is, the number of nonzero entries in $\lambda$.
Using the Markov moves and Equation \eqref{eq:addstrand},
it is easy to show that any normalized Markov trace must satisfy
\begin{equation}
\label{eq:tr}
\tr(b_\lambda) =  \left( \frac{1+q_1 q_2}{q_1 + q_2} \right)^{k-1} \tr(\id_1),
\end{equation}
where $\id_1$ is the identity element of $\Hecke_1$.

We now show that this equation
defines a normalized Markov trace.
Let $b$ be a braid diagram.
Let $v$ be the closed braid diagram obtained by
identifying the top and the bottom of $b$.
Apply the algorithm from the previous section
to write $v$ as a linear combination of $v_\lambda$.
Let $\tr(b)$ be as given by
Equation \eqref{eq:tr} with $\tr(\id_1) = 1$.

Let $v_0$ be the closed braid diagram obtained by
identifying the top and the bottom of $\iota(b)$.
This is obtained from $v$ by adding a disjoint loop.
Now apply the algorithm to write $v_0$
as a linear combination of basis elements $v_\lambda$.
The added loop in $v_0$ remains unchanged throughout the algorithm.
Thus it has the effect of
adding an extra component to each term of
the resulting linear combination of basis elements.
This shows that $\tr$ satisfies Equation \eqref{eq:addstrand}.

Now let $v_+$ be the closed braid diagram obtained by
identifying the top and the bottom of $\sigma_n \iota(b)$.
then $v_+$ is obtained from $v$
by adding an extra ``kink'' in one of the strands.
Now apply the algorithm to write $v_+$
as a linear combination of basis elements $v_\lambda$.
We can assume that we
never choose a basepoint that lies on the added kink.
Then the added kink remains unchanged throughout the algorithm.
It has no effect on the number of components in each term
of the resulting linear combination of basis elements.
This shows that $\tr$ is invariant under the second Markov move.

This completes the proof that $\tr$ is a normalized Markov trace.
Any other normalized Markov trace must be a scalar multiple of $\tr$.
The HOMFLYPT polynomial $P_L(q_1,q_2)$ of a link $L$
is defined to be $\tr(b)$ for any braid $b$ whose closure is $L$.
There are many different definitions of $P_L$ in the literature,
each of which can be obtained from any other by a change of variables.
They are usually specified
by giving the coefficients of the three terms
of the skein relation in Figure \ref{fig:skein}.
As far as I know,
mine is yet another addition
to the collection of possible choices that appear in the literature.

The Jones polynomial $V_L$ is given by
$$V_L(t) = P_L(-t^\frac{1}{2}, t^\frac{3}{2}).$$
If $L$ is a knot,
$V_L$ turns out to involve only integer powers of $t$.
This polynomial was originally defined
as a trace function of the Temperley-Lieb algebra,
which is a certain quotient of the Iwahori-Hecke algebra

A somewhat tangential question is worth mentioning here.
In its most open-ended form, it is as follows.

\begin{qn}
\label{qn:markov}
What are the equivalence classes of braids modulo the moves
\begin{itemize}
\item $ab \leftrightarrow ba$, and
\item $b \leftrightarrow \sigma_n \iota(b)$?
\end{itemize}
\end{qn}

In other words,
what happens if the Markov move
$$b \leftrightarrow \sigma_n^{-1} \iota(b)$$
is omitted?
This question was shown in \cite{OS03}
to be equivalent to the important problem in contact geometry
of classifying transversal links up to transversal isotopy.

The {\em Bennequin number} of a braid $b \in B_n$ is $e-n$,
where $e$ is the sum of the exponents in a word
in the generators $\sigma_i$ representing $b$.
The Bennequin number is invariant under the moves in Question \ref{qn:markov}.
Thus it can be used to show that
there are braids that are related by Markov moves,
but not by the moves in Question \ref{qn:markov}.

Birman and Menasco \cite{BM03} and Etnyre and Honda \cite{EH03}
have independently found pairs of braids
that are related by Markov moves and have the same Bennequin invariant,
but are not related by the moves in Question \ref{qn:markov}.
Their proofs are quite complicated,
and it would be nice to have a new invariant
that could distinguish their pairs of braids.

\section{Representations of $\Hecke_n$}
\label{sec:rep}

A {\em representation} of $\Hecke_n$ is simply a $\Hecke_n$-module.
If $R$ is a field then
an {\em irreducible} representation of $\Hecke_n$ is
a nonzero $\Hecke_n$-module with no nonzero proper submodules.
In \cite{DJ86}, Dipper and James 
gave a complete list of the irreducible representations of $\Hecke_n$.
The aim of this section is to summarize their results.
Our approach comes from the theory of {\em cellular algebras},
as defined in \cite{GL96}.
For convenience we will take
$$\Hecke_n = \Hecke_n(-1,q)$$
from now on.

Let $\lambda$ be a partition of $n$.
The {\em Young subgroup} $\Sym_\lambda$ of $\Sym_n$
is the image of the obvious embedding
$$\Sym_{\lambda_1} \times \dots \times \Sym_{\lambda_k} \to \Sym_n.$$
More precisely,
it is the set of permutations of $\{1,\dots,n\}$ that
fix setwise each set of the form $\{k_i+1,\dots,k_i+\lambda_i\}$
where $k_i = \lambda_1+ \dots + \lambda_{i-1}$.
Let
$$m_\lambda = \sum_{w \in \Sym_\lambda} T_w.$$

Let $\lambda$ be a partition of $n$.
Let $M^\lambda$ be the left-ideal $\Hecke_n m_\lambda$.
We say a partition $\mu$ of $n$ {\em dominates} $\lambda$
if $\sum_{i=1}^j \mu_i \le \sum_{i=1}^j$ for all $j \ge 1$.
Let $I^\lambda$ be the two-sided ideal of $\Hecke_n$
generated by $m_\mu$ for all partitions $\mu$ of $n$ that dominate $\lambda$.
The {\em Specht module} is the quotient
$$S^\lambda = M^\lambda / (M^\lambda \cap I^\lambda).$$

Let $\rad S^\lambda$ be the set of $v \in S^\lambda$
such that $m_\lambda h v = 0$ for all $h \in \Hecke_n$.
Let
$$D^\lambda = S^\lambda / \rad S^\lambda.$$

\begin{thm}
Suppose $R$ is a field.
Then every irreducible representation of $\Hecke_n$
is of the form $D^\lambda$ for some partition $\lambda$ of $n$.
If $\lambda$ and $\mu$ are distinct partitions of $n$
then $D^\lambda$ and $D^\mu$ are either distinct or both zero.
\end{thm}

The definition of $D^\lambda$ can be better motivated
by defining a bilinear form on $S^\lambda$.
Let $\star \co \Hecke_n \to \Hecke_n$ be
the antiautomorphism given by $T_w^* = T_{w^{-1}}$ for all $w \in \Sym_n$.
Note that $m_\lambda^* = m_\lambda$.
The following lemma is due to Murphy \cite{gM92}.

\begin{lem}
\label{lem:mhm}
If $h \in \Hecke_n$ then
$m_\lambda h m_\lambda = r m_\lambda$ modulo $I^\lambda$,
for some $r \in R$.
\end{lem}

Thus we can define a bilinear form
$$\langle \cdot , \cdot \rangle \co S^\lambda \times S^\lambda \to R$$
by
$$(h_1 m_\lambda)^* (h_2 m_\lambda)
         = \langle h_1 m_\lambda, h_2 m_\lambda \rangle m_\lambda.$$
Then $\rad S^\lambda$ is the set of $y \in S^\lambda$
such that $\langle x,y \rangle = 0$ for all $x \in \Hecke_n$.

Dipper and James also determined
which values of $\lambda$ give a nonzero $D^\lambda$.
Let $e$ be the smallest positive integer such that $1+q+ \dots + q^{e-1} = 0$,
or infinity if there is no such integer.

\begin{thm}
$D^\lambda \neq 0$ if and only if
$\lambda_i - \lambda_{i+1} < e$ for all $i \ge 1$.
\end{thm}

Thus the work of Dipper and James completely characterizes
the irreducible representations of $\Hecke_n$.
However, understanding these irreducible representations
remains an active area of research to this day.
An example of a major open-ended question in the area is the following.

\begin{qn}
What can be said about the dimensions of $D^\lambda$?
\end{qn}

\section{The future}
\label{sec:future}

My hope for the future is that
the definition of $B_n$ as a mapping class group
will provide solutions to problems related to
the representation theory of the Iwahori-Hecke algebra.
One reason for optimism is the mechanism described in \cite{sB04}
to obtain representations of the $\Hecke_n$
from the induced action of $B_n$
on homology modules of configuration spaces in $D_n$.
There I conjectured that
all irreducible representations $D^\lambda$
can be obtained in this way.
The inner product we defined on $D^\lambda$
would presumably correspond to
the intersection form on homology.

Another direction for future research
is to look at other quotient algebras of $RB_n$.
After the Iwahori-Hecke algebra,
the next obvious candidate
is the Birman-Wenzl-Murakami algebra.
A somewhat non-standard presentation of this algebra
is as follows.

Let $X$ be the following element of the braid group algebra $RB_n$.
$$X = q \bar \sigma_1 + 1 - q - \sigma_1.$$
The Birman-Wenzl-Murakami algebra
is the quotient of $RB_n$
by the following relations.
\begin{itemize}
\item $(q^2 \sigma_1^{-1} \sigma_2^{-1} - \sigma_1 \sigma_2)X = 0$,
\item $(q \sigma_2^{-1} + 1 - q - \sigma_2) X =
       (q \sigma_1^{-1} \sigma_2^{-1} - \sigma_1 \sigma_2) X$,
\item $\sigma_1 X = t X$.
\end{itemize}

This algebra has an interesting history.
After the discovery of the Jones polynomial,
Kauffman \cite{lK90} discovered a new knot invariant
which he defined directly using the knot or link diagram
and a skein relation.
The Birman-Wenzl-Murakami algebra was then constructed
in \cite{BW89}, and independently in \cite{jM87},
so as to give the Kauffman polynomial via a trace function.
Thus the history of the Birman-Wenzl-Murakami algebra
traces the history of the Jones polynomial in reverse.

\begin{qn}
How much of this paper
can be generalized to the Birman-Wenzl-Murakami algebra?
\end{qn}

In this direction, John Enyang \cite{jE04} has shown that
the Birman-Murakami-Wenzl algebra is a cellular algebra,
and used this to give a definition of its irreducible representations
similar to the approach in Section \ref{sec:rep}.

Next we would like to generalize \cite{sB04}
to the Birman-Wenzl-Murakami algebra.

\begin{qn}
Is there a homological definition
of representations of the Birman-Wenzl-Murakami algebra?
\end{qn}

I believe the answer to this is yes.
Furthermore,
the homological construction suggests a new algebra $Z_n$,
which would further generalize
the Iwahori-Hecke and Birman-Wenzl-Murakami algebras.
I will conclude this paper with a definition of $Z_n$
and some related open questions.
I hope these might be amenable to some combinatorial computations,
even without the homological motivation,
which is currently unclear and unpublished.

We use the notation
$\sigma_{i_1 \dots i_k}$ as shorthand for $\sigma_{i_1} \dots \sigma_{i_k}$,
and $\bar \sigma_{i_1 \dots i_k}$ for $\sigma_{i_1 \dots i_k}^{-1}$.
Define the following elements of $RB_n$.
\begin{eqnarray*}
X_2 &=& q \bar \sigma_1 + 1 - q - \sigma_1      \\
X_3 &=& (q^2 \bar \sigma_{21} - \sigma_{12})X_2 \\
X_4 &=& (q^3 \bar \sigma_{321} - \sigma_{123})X_3 \\
    &\vdots& \\
X_n &=& (q^{n-1} \bar \sigma_{(n-1)\dots 1} - \sigma_{1 \dots (n-1)}) X_{n-1}
\end{eqnarray*}
Then $Z_n$ is the algebra $RB_n$ modulo the following relations.
\begin{eqnarray*}
(q \bar \sigma_2 + 1 - q - \sigma_2) X_2 &=&
(q \bar \sigma_{21} - \sigma_{12}) X_2,
\\
(q^2 \bar \sigma_{32} - \sigma_{23}) X_3 &=&
(q^2 \bar \sigma_{321} - \sigma_{123}) X_3,
\\
(q^3 \bar \sigma_{432} - \sigma_{234}) X_4 &=&
(q^3 \bar \sigma_{4321} - \sigma_{1234}) X_4,
\\
&\vdots&
\\
(q^{n-1} \bar \sigma_{n \dots 2} - \sigma_{2 \dots n}) X_n &=&
(q^{n-1} \bar \sigma_{n \dots 1} - \sigma_{1 \dots n}) X_n.
\end{eqnarray*}

Note that $\Hecke_n(1,-q)$ is the quotient of $Z_n$
by the relation $X_2 = 0$.
Also the Birman-Wenzl-Murakami algebra
is the quotient of $Z_n$
by the relations
$X_3 = 0$ and $\sigma_1 X_2 = t X_2$.
The following basically asks if $Z_n$
is bigger than the Birman-Wenzl-Murakami algebra.

\begin{qn}
Does $X_3$ equal $0$ in $Z_n$?
\end{qn}

Presumably some extra relations should be added to $Z_n$,
such as $\sigma_1 X_2 = t X_2$,
or something more general.

\begin{qn}
What extra relations should be added to $Z_n$ to make it finite-dimensional?
\end{qn}

\begin{qn}
How much of this paper can be generalized to $Z_n$?
\end{qn}

It might be easier to first study these questions
for the quotient of $Z_n$ by the relation $X_4 = 0$.


\newcommand{\etalchar}[1]{$^{#1}$}
\providecommand{\bysame}{\leavevmode\hbox to3em{\hrulefill}\thinspace}
\providecommand{\MR}{\relax\ifhmode\unskip\space\fi MR }
\providecommand{\MRhref}[2]{%
  \href{http://www.ams.org/mathscinet-getitem?mr=#1}{#2}
}
\providecommand{\href}[2]{#2}


\end{document}